*To my high school math teachers Carl Moore and Dennis Fatland*

# PROJECTIVE LINE REVISITED

NICHOLAS PHAT NGUYEN

**Abstract.** This article provides a new perspective on the geometry of a projective line, which helps clarify and illuminate some classical results about projective plane. As part of the same train of ideas, the article also provides a proof of the nine-point circle theorem valid for any affine plane over any field of characteristic different from 2.

**1. INTRODUCTION.** The projective line is such a humble and well-known mathematical object that there seems to be hardly anything interesting to talk about. However, in this note, we want to show the readers a different way to think about the geometry of a projective line that the readers will not see in any projective geometry textbook. This new perspective will greatly clarify a number of classical results about the geometry of a projective plane, and will also help us outline a proof of the nine-point circle theorem that applies to any affine plane over a field of characteristic $\neq 2$.

Let $K$ be any field. The projective line $K \cup \{\infty\}$ can be described as the affine line $K$ extended by adjoining a point at infinity denoted by the symbol $\infty$. It can also be regarded as $\mathbf{P}(K^2)$ = the set of lines through the origin in the affine plane $K^2$, where each element $x$ of $K$ is identified with the line passing through the point $(x, 1)$, and $\infty$ is identified with the horizontal line consisting of all the points $(u, 0)$ ($u$ running through $K$) in $K^2$.

A projective transformation of a projective line is a transformation of $K \cup \{\infty\}$ induced by an invertible linear transformation of $K^2$ (viewed as a vector space of dimension 2 over $K$). All such projective transformations can be described as the fractional linear transformations of $K \cup \{\infty\}$ described by:

$$x \to x' = \frac{ax+b}{cx+d} \text{ where } a, b, c, d \text{ are numbers in } K \text{ such that } ad - bc \neq 0,$$

A major focus of the classical geometry of projective line is the study of involutions, defined as projective transformations of order two. These involutions are precisely the fractional linear transformations where $a + d = 0$.[2] In this note, we will outline a different way to view these involutions when the field $K$ has characteristic $\neq 2$. Such a new perspective can take us a long way, as we hope to show.

**2. INVOLUTIONS**. As explained below, we can think of involutions on a projective line as the 1-dimensional versions of the familiar inversions and reflections in Euclidean plane geometry.

Recall that in the extended Euclidean plane $\mathbf{R}^2 \cup \{\infty\}$, we can define reflections (across a line) and inversions (across a circle) as follows.

<u>Reflection</u>: Consider the line defined by the equation $L(\mathbf{x}) = \mathbf{b}.\mathbf{x} + c = 0$ where $\mathbf{b}$ is a nonzero vector and $\mathbf{b}.\mathbf{x}$ is the standard dot product of two vectors. The corresponding affine reflection is the transformation of $\mathbf{R}^2 \cup \{\infty\}$ that sends $\infty$ to itself, and any finite point $\mathbf{x}$ to a finite point $\mathbf{x'} = \mathbf{x} - 2\mathbf{b}\frac{L(\mathbf{x})}{\mathbf{b}.\mathbf{b}}$. Note that $\mathbf{b}.\mathbf{b}$ is a nonzero number because $\mathbf{b}$ is a nonzero vector and the standard dot product in $\mathbf{R}^2$ is anisotropic.

<u>Inversion</u>: Any circle in $\mathbf{R}^2$ can be defined by an equation $a(\mathbf{x} - \mathbf{b}).(\mathbf{x} - \mathbf{b}) - as = 0$, where $a$ is a nonzero real number. The center of the circle is the point represented by the vector $\mathbf{b}$, and the radius of the circle is zero when $s$ is zero. Take any circle with nonzero radius, i.e. $s$ is nonzero. For such a circle, we can define an inversion mapping $\mathbf{R}^2 \cup \{\infty\}$ to itself as follows:

- $\mathbf{b} \leftrightarrow \infty,$

---

[2] For a projective transformation to be an involution, it is necessary and sufficient that any 2x2 matrix T corresponding to such an involution satisfies a minimal equation $T^2 - s$ with $s \neq 0$. At the same time, we know from the Hamilton-Cayley theorem that $T^2 - tr(T)T + det(T) = 0$. So the trace of such a matrix must be zero.

- for any finite point $x \neq b$, $x \rightarrow$ the finite point $x'$ collinear with $b$ and $x$, and such that $(x - b).(x' - b) = s$.

The collinear condition and dot product equation are equivalent to $x' - b = \dfrac{s(x-b)}{(x-b).(x-b)}$.

Note that $(x - b).(x - b)$ is a non-zero number because $(x - b) \neq 0$. Incidentally, this shows that inversion is a well-defined transformation regardless of whether or not there is any point $x$ in the plane satisfying the circle equation $a(x - b).(x - b) - as = 0$.

Any non-zero scalar multiple of a line or circle equation defines the same transformation, and any reflection or inversion so defined is a transformation of $\mathbb{R}^2 \cup \{\infty\}$ that is its own inverse. If $x$ and $x'$ are mapped to each other by such a transformation, we will call these points inverse points or conjugate points relative to the defining line or circle.

For a projective line over any field $K$ of characteristic $\neq 2$, we can define reflections and inversions in the same way.

<u>Reflection on a projective line</u>: The analog of a line in this setting is the linear expression $Bx + C$ in one variable x, where $B$ is a nonzero number in $K$. The equation $Bx + C = 0$ has exactly one zero in $K$, namely the number $(-C/B)$. The reflection defined by the line $Bx + C$ is the transformation that sends $\infty$ to itself, and any finite point $x$ to the finite point $x^* = x - 2B\dfrac{Bx+C}{B.B} = -x - 2(C/B)$. (In this context, the dot product notation for two numbers just means ordinary multiplication.) Such a transformation is a combination of the symmetry $x \rightarrow -x$ and the translation $x \rightarrow x - 2(C/B)$. The fixed points of such a reflection are the finite point $-C/B$ and the point $\infty$ at infinity.

<u>Inversion on a projective line</u>: The analog of a circle in this setting is the quadratic expression $Ax^2 + Bx + C$, where the leading coefficient $A$ is a nonzero number. We can rewrite the expression as $A(x + B/2A)^2 - As$ where $s = (B^2 - 4AC)/4A^2$. When the discriminant $(B^2 - 4AC)$ of this quadratic expression vanishes, we have the analog of a zero circle (circle with zero radius). When the discriminant $(B^2 - 4AC)$ is nonzero, we can define an inversion of the projective line as the transformation $x \rightarrow x^*$ that exchanges the points $\infty$ and − *B/2A*, and otherwise *(x\* + B/2A)(x + B/2A) = (B² − 4AC)/4A²*. Note that this transformation is well-defined when the right-hand side is nonzero, regardless of whether the quadratic polynomial *Ax² + Bx + C* has any root in *K*. If that quadratic polynomial has two distinct roots in *K*, then these two roots are fixed points of the transformation.

From the above definitions, it is clear that reflections and inversions are involutions. Moreover, any involution is either a reflection or an inversion. Indeed, an involution on a projective line is just a fractional linear transformation $x \to x^* = \frac{ax+b}{cx+d}$ where *a + d = 0*. If *c = 0*, then we have a reflection $x \to x^* = -x + \frac{b}{d}$. If *c ≠ 0*, then we have an inversion defined by the equation *c(x\* + d/c)(x + d/c) = (d² + bc)/c*, i.e., the inversion across the circle *cx² + 2dx − b*.

For the rest of our discussion, we will always assume that the underlying field K has characteristic ≠ 2.

**3. THE SPACE OF INVOLUTIONS.** Seen as a reflection or an inversion, each involution on a projective line is represented (uniquely up to a scalar factor) by a linear or quadratic expression. Consider the set *E* of all polynomials *p(X)* of degree ≤ 2, with coefficients in the field *K*. The set *E* is naturally a *K*-vector space of dimension 3. We will refer to a nonzero polynomial *p* as a 2-cycle, 1-cycle, or 0-cycle depending on whether the degree of *p* is 2, 1, or 0. For convenience, we will write each element *p* of *E* in the same form *p(X) = aX² + bX + c*, with the understanding that each coefficient *a*, *b*, and *c* could be zero.

We can endow the vector space *E* with a symmetric bilinear form < _ , _ > as follows. Given cycles *p = aX² + bX + c* and *p\* = a\*X² + b\*X + c\**, we define *<p,p\*>* to be *bb\* − 2ac\* − 2a\*c*.

This scalar product is clearly symmetric and bilinear. Moreover, it is nondegenerate, because it is plainly isomorphic to the sum of *K* (represented by the middle coefficient, with ordinary multiplication) and an Artinian plane (also known as a hyperbolic plane). We will refer to this fundamental scalar product on *E* as the cycle pairing or cycle product.

The vector space $E$ can be naturally identified with the vector space $Q$ of all symmetric bilinear forms on $K^2$. Specifically, an element $p$ of $E$ can be thought of as a function from $K$ to $K$ given by $p(x) = q((x, 1), (x, 1))$ where $q$ is a symmetric bilinear form on $K^2$.[3] If $p = Ax^2 + Bx + C$, then the matrix for the corresponding symmetric bilinear form $q$ (relative to the standard basis of $K^2$) has entries $A$ and $C$ in the main diagonal, and entries $B/2$ in the cross diagonal. For convenience, we will often write a symmetric bilinear form as a homogenous polynomial of degree 2, so that the bilinear form $q$ corresponding to the cycle $p$ above is $Ax^2 + Bxy + Cy^2$.

The cycle pairing defined on $E$ can be carried over to a regular scalar product on the vector space $Q$. If we express a symmetric bilinear form in $Q$ as a 2 x 2 matrix, then the norm of such a matrix under this cycle product is simply $-4$ times the determinant. Specifically, for a cycle $p = Ax^2 + Bx + C$ in $E$, its norm $<p,p> = B^2 - 4AC$ is equal to $-4$ times the determinant of the corresponding bilinear form $q$.

For computational purposes, it is often easier to work with cycles. However, because symmetric bilinear forms have an intrinsic meaning independent of coordinates, it can be helpful sometimes to think of cycles in terms of symmetric bilinear forms. For example, if we have a change in coordinates, what will happen to the cycle pairing? If we think of cycles as symmetric bilinear forms, the question has a straight-forward conceptual answer, as described below.

If we change the coordinates for $\mathbf{P}(K^2)$ by means of a general linear transformation $x = S(x')$ of $K^2$, then the matrix $M$ of a symmetric bilinear form in the old coordinate $x$ will become ${}^tSMS$ (where ${}^tS$ is the transpose of $S$) in the new coordinates $x'$. The norm $<M, M>$ of $M$ under the cycle pairing is equal to $-4\det(M)$. With this change in coordinates, the norm of $M$ becomes $<{}^tSMS, {}^tSMS> = -4\det({}^tSMS) = -4\det(S)^2\det(M)$. So the simple linear transformation $q \mapsto (\det S)q$ gives us an isometry between the vector space $Q$ with the cycle pairing in the new coordinate $x'$ and the vector space $Q$ with the cycle pairing in the old coordinate $x$. Accordingly, orthogonal properties in the space $Q$ under the cycle pairing are independent of any coordinate chosen for the parametrization of the projective line.

---

[3] The elements $p$ of $E$ are defined as polynomials of degree 2 or less, but because the field $K$ has 3 or more elements, such a polynomial $p$ can be identified with a polynomial function from $K$ to $K$.

Any 1-cycle $b(X - u)$ seems to have just one zero point, namely $u$. However, if we think of that 1-cycle as equivalent to the bilinear form $b(XY - uY^2)$, then we have two linearly independent isotropic vectors $(u, 1)$ (corresponding to the point $u$) and $(x, 0)$ corresponding to the point at infinity). Because of this fact, we will regard $\infty$ as the second zero point of any 1-cycle. Similarly, we regard $\infty$ as the zero point of any 0-cycle $C$ ($C$ a nonzero constant) in light of the fact that the corresponding bilinear form $CY^2$ has any $(x, 0)$ as an isotropic vector.

With this convention, each involution is identified (up to a scalar factor) with a nonisotropic cycle whose zero points (if any) are the fixed points of the involution. For example, the zero points of the 1-cycle $(X - u)$ are the points $u$ and $\infty$, which are exactly the fixed points of the reflection defined by the 1-cycle $(X - u)$. For the 2-cycle $(X - u)(X - w)$, its zero points are the distinct points $u$ and $w$, which are exactly the fixed points of the inversion defined by the 2-cycle $(X - u)(X - w)$.

The isotropic elements of $E$ (with respect to the cycle pairing defined above) are the 0-cycles and the 2-cycles with zero discriminant. The nonisotropic elements of $E$ are the 1-cycles and the 2-cycles with nonzero discriminant, precisely the elements for which we can define reflections and inversions. In $Q$, the corresponding isotropic elements are the degenerate bilinear forms, and the nonisotropic elements are the nondegenerate bilinear forms. Because each involution on a projective line is represented by a nonisotropic cycle in $E$ uniquely up to a scalar factor, we can regard all the nonisotropic elements in the projective space $\mathbf{P}(E)$ or $\mathbf{P}(Q)$ as the space of all involutions for the projective line $K \cup \{\infty\} = \mathbf{P}(K^2)$.

We say that a 2-dimensional subspace of $E$ or $Q$ is regular if the scalar product induced by the cycle pairing on that subspace is nondegenerate. The corresponding projective line (also called pencil) in the projective subspace $\mathbf{P}(E)$ or $\mathbf{P}(Q)$ is then also said to be regular.

Because the cycle pairing is nondegenerate, there is a natural bijective correspondence between nonisotropic elements of $P(E)$ or $P(Q)$ and regular pencils of cycles or bilinear forms. Specifically, the orthogonal complement of any nonisotropic element is a regular pencil, and vice versa.

**Proposition 1:** *There is a natural bijection between involutions of a projective line and regular pencils of bilinear forms on that projective line. Each pair of conjugate points in an involution are the zero points of a symmetric bilinear form in the regular pencil corresponding to that involution.*

*Proof.* Recall that each involution is either a reflection or an inversion. These involutions are parametrized by the set of nonisotropic elements in the projective space $P(E)$ or $P(Q)$. Moreover, each such nonisotropic element corresponds exactly to a regular pencil if we look at the orthogonal complement. So involutions on a projective line correspond bijectively to regular pencils of bilinear forms.

Consider first the case of a reflection defined by the 1-cycle $(X - c)$. The fixed points of this reflection are the finite point $c$ and the point at infinity $\infty$. The finite point $c$ is the zero point of the 2-cycle $(X - c)^2 = X^2 - 2cX + c^2$. The cycle pairing of that 2-cycle with the 1-cycle $(X - c)$ is simply $-2c + 2c = 0$. Similarly, the point at infinity is the zero point of any 0-cycle, which is easily seen to be orthogonal to any 1-cycle. Accordingly, the proposition is certainly true for the fixed points.

If $u$ and $w$ are two distinct points that are conjugate under this reflection, then they are finite points such that $w = -u + 2c$ or $u + w = 2c$. The points $u$ and $w$ are the zero points of the 2-cycle $(X - u)(X - w)$. We need to show that the pairing of the 1-cycle $(X - c)$ and the 2-cycle $(X - u)(X - w)$ is zero. But under the formula for the cycle pairing product, their cycle pairing is simply just $-(u + w) + 2c = 0$.

Now consider the case when the involution in question is an inversion defined by a 2-cycle $(X^2 - 2bX + c)$ with nonzero discriminant $4b^2 - 4c$. The conjugate points $\infty$ and $b$ are the zero points of the 1-cycle $(X - b)$. The 1-cycle $(X - b)$ and the 2-cycle $(X^2 - 2bX + c)$ have the pairing $-2b + 2b = 0$ and so are clearly orthogonal.

Other pairs of conjugate points $u$ and $w$ are related by the equation $(u - b)(w - b) = b^2 - c$. That can be written as $uw - b(u + w) + c = 0$. We want to show that the pairing of the 2-cycle $(X^2 - 2bX + c)$ and the 2-cycle $(X - u)(X - w)$ is zero. But that pairing is $2b(u + w) - 2uw - 2c = 0$ in light of the relationship between $u$ and $w$. ∎

In the above correspondence, the fixed points of an involution correspond (up to a scalar factor) to degenerate bilinear forms. Because a reflection or an inversion on a projective line will either have no fixed point or exactly two fixed points, a regular pencil of bilinear forms will either have no degenerate element or exactly two degenerate elements. This is a geometrical interpretation of the well-known fundamental result that a regular symmetric bilinear space of dimension 2 is either anisotropic or an Artinian plane (with two linearly independent isotropic vectors).

Another consequence of Proposition 1 is that two involutions with the same two fixed points must be the same transformation. That is because the bilinear forms corresponding to these fixed points will generate the same pencil of bilinear forms, and therefore we must have the same involution.

**4. THE INVOLUTION THEOREM OF DESARGUES.** This equivalence between involutions and regular pencil of symmetric bilinear forms gives us a precise and general condition for the involution theorem of Desargues. Girard Desargues, one of the founders of projective geometry, discovered the following remarkable theorem, which can be stated roughly as follows: *A pencil of conics in a projective plane will generally intersect a line in pairs of points that are conjugate under an involution.*

We say "generally" because this theorem is true for most but not all configurations. If the pencil in question is the set of all conics passing through 4 points in general position (meaning in this case that no 3 of them are collinear), then the most commonly stated condition is that the line does not pass through any of the 4 given points.

Proposition 1 gives a precise condition for the Desargues involution theorem to be true in the case of any general pencil of conics. A pencil of conics in a projective plane corresponds to a two-dimensional space of symmetric bilinear forms. For that pencil of conics to intersect a given projective line in conjugate pair of points under an involution, Proposition 1 implies that the bilinear forms in that pencil when restricted to the line must be a regular two-dimensional space with respect to a cycle pairing on that projective line.[4]

For example, consider the projective plane with homogenous coordinates $(u, v, w)$ and assume that the projective line in question is the line given by $w = 0$. Each conic $Au^2 + Buv + Cv^2 +$ (terms with variable $w$) becomes the bilinear form $Au^2 + Buv + Cv^2$ when restricted to the line $w = 0$. The cycle pairing on the line $w = 0$ is the following. If $f = Au^2 + Buv + Cv^2$ and $g = au^2 + buv + cv^2$, then we have $<f, g> = Bb - 2Ac - 2Ca$. The Desargues involution theorem holds if and only if the pencil of bilinear forms as restricted to the given line is regular with respect to the above scalar product.

In general, it is straight-forward to check whether a bilinear space of dimension 2 is regular under a given symmetric pairing. We can, for example, just write down the matrix of that pairing relative to a suitable basis of the space and determine if the matrix has non-zero determinant. In the particular case of the space $E$ or $B$, we have another geometric criterion.

**Proposition 2:** *A 2-dimensional subspace of E is regular with respect to the cycle pairing if and only if there is no common zero point for all the cycles in that subspace, or equivalently, if the subspace can be generated by two cycles with no common zero point.*

*Proof.* Let $L$ be a 2-dimensional subspace of $E$. Recall that $E$ has dimension 3 and is regular under the cycle pairing. Because $E$ is regular, the subspace $M$ of $E$ orthogonal to $L$ is therefore a 1-dimensional subspace, say generated by a cycle $h$. It follows that $L$ is a regular subspace of $E$ if and only if $h$ is outside L. If $h$ also belongs to $L$, then it is an isotropic cycle orthogonal to all cycles in $L$. In that case, $h$ is either a 2-cycle centered at a point $u$ of $K$ that is a common zero for all cycles in $L$, or $h$ is a 0-cycle whose zero point (the point at infinity $\infty$) is a common zero for all the cycles in *L* (in which case *L* must necessarily be the 2-dimensional subspace of all 1-cycles and 0-cycles in *E*). ∎

---

[4] Recall that such cycle pairings under different coordinates are all isomorphic.

Based on Proposition 2, a family of all conics passing through four points in general position in a projective plane would induce a regular 2-dimensional subspace of cycles on a projective line (and therefore an involution on that line under proposition 1) if and only if these conics have no common zero on that projective line, i.e. if and only if the projective line does not pass through any of the four base points of that family.

5. **THE ELEVEN-POINT CONIC.** A nondegenerate symmetric bilinear form on the two-dimensional vector space $K^2$ gives us a natural involution of the projective line $P(K^2)$ which maps a point on the line to its polar conjugate. This is well-defined because the bilinear form is nondegenerate. We will refer to this involution as the polar involution defined by a nondegenerate bilinear form.

The orthogonal complement of such a nondegenerate bilinear form is a regular pencil of bilinear forms, and therefore gives us an involution on the projective line under correspondence of Proposition 1. We will refer to this involution as the Desargue involution defined by a nondegenerate bilinear form.

**Proposition 3**. *The Desargue involution and the polar involution defined by a nondegenerate symmetric bilinear form are the same transformation.*

*Proof.* Let the nondegenerate symmetric bilinear form be $q((x,y),(x',y')) = Axx' + Byx' + Bxy' + Cyy'$, and let $(u, v)$ be the homogeneous coordinates of the projective line in question. The polar involution $(u,v) \mapsto (u^*, v^*)$ is defined by the equation:

$$q((u,v), (u^*, v^*)) = 0 = Auu^* + Bvu^* + Buv^* + Cvv^*$$

For the Desargues involution, Proposition 1 tells us that the two conjugate points $(u,v)$ and $(u^*, v^*)$ are isotropic points of a bilinear form *h* which is cycle-orthogonal to *q*.

Now note that the two points $(u, v)$ and $(u^*, v^*)$ are the zero or isotropic points of $h$ and also of the bilinear form $(vX - uY)(v^*X - u^*Y)$. If two bilinear forms of dimension 2 have the same two isotropic points, then they must be proportional. Indeed, relative to the basis consisting of those two isotropic vectors the 2 x 2 matrices of these two bilinear forms both have zeros in the diagonal and a non-zero number in the cross diagonal.

Accordingly, the bilinear form $(vX - uY)(v^*X - u^*Y)$ must also be cycle-orthogonal to $q$. Writing out $(vX - uY)(v^*X - u^*Y) = vv^*X^2 + (-vu^* - uv^*)XY + uu^*Y^2$, the cycle-orthogonal relationship means that we have the equation $2B(-vu^* - uv^*) - 2Auu^* - 2Cvv^* = 0$. This equation is the same as the equation for the polar involution up to a factor of $-2$. This means the Desargue involution and the polar involution are the same transformation. ∎

Proposition 3 allows us to gain some more insight into the following remarkable conic. Consider the pencil $L$ of conics passing through four points in general position in a projective plane. Suppose that this pencil when restricted to a line $D$ gives us a regular pencil of bilinear forms on that projective line, so that the Desargues involution theorem holds for that pencil $L$ and the line $D$. For each non-degenerate conic in the pencil $L$, consider the pole of $D$ relative to that conic. All such poles as the conics range over the pencil $L$ constitute a conic $\mathcal{E}$ that has some remarkable properties. See, e.g., [3] at section 79.1. In particular, that conic $\mathcal{E}$ passes through potentially up to eleven points that are defined by the configuration of $D$ and the base points of the pencil $L$. Accordingly, it is known as the eleven-point conic.

By hypothesis the conics of the pencil $L$ intersect the line $D$ (when they do) in pairs of conjugate points under a Desargue involution. In addition, we also have an involution of the line $D$ induced by the conic $\mathcal{E}$.

**Proposition 4**. *The eleven-point conic $\mathcal{E}$ induces a nondegenerate bilinear form on the line $D$, and therefore gives us an involution defined by polarity with respect to $\mathcal{E}$. The polar involution induced by the eleven-point conic $\mathcal{E}$ on the line $D$ is the same as the Desargue involution induced by the pencil $L$ on the line $D$.*

*Proof.*   Note that any isotropic point of the conic $\mathcal{E}$ on the line $D$ is exactly a tangent point of $D$ with a conic in the pencil $L$. Each such tangent point corresponds to a degenerate bilinear form in the pencil $L$ when restricted to $D$ (since the tangent point shows that the corresponding bilinear form has a nonzero radical). Because the pencil $L$ when restricted to $D$ is a regular pencil of bilinear forms, we either have no such tangent point, or exactly two tangent points. Accordingly, the conic $\mathcal{E}$ when restricted to the line $D$ will either have no isotropic point (when no conic in the pencil $L$ is tangent to $D$), or exactly two isotropic points (when two different conics in the pencil $L$ are tangent to $D$). In other words, the conic $\mathcal{E}$ when restricted to the line $D$ also gives us a nondegenerate bilinear form.

If the Desargue involution induced by the pencil $L$ on the line $D$ has two distinct fixed points, then these two fixed points are also zero points of the conic $\mathcal{E}$ because they are necessarily the tangent points of the line $D$ with two conics in the pencil $L$. These fixed points are therefore the intersection points of the line $D$ with the conic $\mathcal{E}$, and consequently are also invariant under the polar involution on $D$ induced by the conic $\mathcal{E}$. But two involutions of a projective line with the same two fixed points must be the same.

If the Desargue involution induced by the pencil $L$ has no fixed point, then we look at the same configuration and equations in the algebraic closure of the base field $K$. In that algebraic closure, any regular symmetric bilinear space of dimension 2 must have two linearly independent isotropic vectors, and therefore the involution induced by the pencil must have two fixed points. Consequently, by extending the base field to its algebraic closure, we see that the two involutions are the same. But that can only be the case if they are already the same transformations over the base field $K$. ∎

We know that the Desargue involution on the line $D$ is either a reflection or an inversion defined by a nondegenerate symmetric bilinear form $q$ that is cycle orthogonal to the pencil $L$ when restricted to $D$. According to Proposition 3, such a Desargue involution is the same as the polar involution defined by $q$. Proposition 4 tells us that the polar involution defined by the eleven-point conic $\mathcal{E}$ on $D$ is the same as the Desargue involution. That means when restricted to $D$, the conic $\mathcal{E}$ will give us the same bilinear form $q$ up to a scalar factor.

In other words, the conic $\mathcal{E}$ is cycle orthogonal to all conics in the pencil $L$ when restricted to the projective line $D$.

For example, consider a projective plane with homogeneous coordinates $(u, v, w)$ and suppose that the line $D$ is the line at infinity given by the equation $w = 0$. If the conic $\mathcal{E}$ has equation $Au^2 + Buv + Cv^2 +$ (terms with variable $w$), and if $r = au^2 + buv + cv^2 +$ (terms with variable $w$) is any conic in the pencil $L$, then such orthogonal relationship means that we must have $Bb - 2Ac - 2Ca = 0$.

Assume that the bilinear forms in such a pencil have the form $au^2 + buv - av^2 +$ (terms with variable $w$). That is to say, the coefficients of the terms $u^2$ and $v^2$ have opposite signs while the coefficient of the term $uv$ ranges over all values in the coefficient field $K$. In that case, the above orthogonal relationship means that the equation for $\mathcal{E}$ must have the form $Au^2 + Av^2 +$ (terms with variable $w$). In other words, $\mathcal{E}$ must be a circle in this case.

This particular situation happens for the following configuration. Let $M$, $N$, and $P$ be 3 points in the affine plan $(u, v, 1)$ that are not collinear. We have the standard dot product $(u, v).(u', v') = uu' + vv'$ in the vector space $K^2$. This standard dot product is non-degenerate, although it is not anisotropic in general because $K$ is an arbitrary field of characteristic $\neq 2$.

By reference to the above standard dot product, we can define orthogonal lines in the affine plane $(w = 1)$, and through each vertex of the triangle $MNP$ there is a unique line orthogonal to the opposite side called the altitude line. It follows as an exercise in linear algebra that all three altitude lines are concurrent, i.e., they pass through a common point $T$ called the orthocenter of the triangle $MNP$.

We will assume that the orthocenter $T$ does not lie on any of the sides of the triangle $MNP$, so that the 4 points $M$, $N$, $P$, $T$ are in general position and form a frame for the projective plane. In that general case, there is a pencil of conics through the four points $M$, $N$, $P$ and $T$ as base points. That pencil can be generated by linear combinations of the following two conics:

(line equation for $MT$) x (line equation for $NP$), and

(line equation for *NT*) x (line equation for *MP*)

For any two line equations that are orthogonal (relative to the standard dot product on *u* and *v*), their product will be an expression of the form $au^2 + buv - av^2 +$ (terms with variable *w*). Consequently, all of the conics in the pencil will have expressions of the same form.

Because the line at infinity does not pass through any of the four base points in general position, this pencil of conics induces a Desargues involution on the line at infinity. It follows from our earlier analysis that the eleven-point conic $\mathcal{E}$ relative to such a pencil must be a circle. From projective geometry, we know that this circle passes through the following nine points determined by the configuration *M, N, P, T* in the affine plane ($w = 1$) and the line at infinity ($w = 0$): namely the midpoints of the six sides of the configuration (*MN, PT, MP, NT, MT, NP*),[5] and the intersections of the three pairs of lines determined by the configuration (*MN* and *PT*, *MP* and *NT*, *MT* and *NP*). In classical Euclidean geometry, this nine-point circle is known as the Feuerbach circle.

What we have shown is the following generalization of the Feuerbach nine-point circle from classical Euclidean geometry:

**Proposition 5**. *In any affine plane over any field of characteristic $\neq 2$, we have a nine-point circle associated with any triangle whose orthocenter is not collinear with any of the sides, similar to the case of the classical Euclidean plane.* ∎

**REFERENCES**

1. Scharlau, W. (1985). *Quadratic and Hermitian Forms*. Grundlehren der Mathematischen Wissenschaften, vol. 270. Berlin, Heidelberg, New York and Tokyo: Springer-Verlag.

---

[5] The midpoint of a segment *MN* in an affine plane is simply the point defined by barycentric coordinates $(1/2)M + (1/2)N$. The points *M, N*, the midpoint of *MN*, and the intersection of the line through *MN* and the line at infinity form a harmonic range with cross ratio -1.

NICHOLAS PHAT NGUYEN

12015 12th Dr SE, Everett, WA 98208, U.S.A

E-mail: nicholas.pn@gmail.com